\documentclass[12pt]{amsart}
\usepackage[latin1]{inputenc}
\usepackage{amsfonts,amsxtra,amssymb,amsthm,amsmath,amscd,mathrsfs, epsfig,url}
\def\leq{\leqslant}

\def\geq{\geqslant}

\theoremstyle{plain}

\theoremstyle{remark}

\theoremstyle{definition}

\numberwithin{equation}{section}

\begin{document}


\title[ Almost cubes in short intervals]
{Almost $k$-th powers in short intervals}

\author{  Wenguang Zhai}

\address{
Department  of Mathematics
\\
China University of Mining and Technology
\\
Beijing 100083
\\
China} \email{zhaiwg@hotmail.com}

\date{\today}

\begin{abstract} Let $k\geq 3$ be a fixed integer and $x$ be a large real number.
Let $0<c_0<\theta\leq 1$ be a fixed real number for some $c_0>0.$
In this paper we show that  there exists a constant $1/2\leq \delta_k(\theta)<1$ such that the interval $[x, x + x^{\delta_k(\theta) +\varepsilon} ]$ contains an integer of the form $n_1n_2 \cdots n_k$ such that
  $|n_j-n^{1/k}|\ll n^{\theta/k} \ (j=1,2, \cdots, k)$. Especially  we have
  $\delta_3(1)=\delta_4(1)=1/2,$ which improves previous results of Chan.
\end{abstract}

\subjclass[2000]{11N25,11L07}
\keywords{almost $k$-th power, exponential sum, prime numbers, short interval}

\thanks{This work is  partially supported  by
the National Natural Science Foundations of China(Grant Nos. 12471009,  12301006)
and  partially supported  by  Beijing Natural Science Foundation (Grant No. 1242003).}

\maketitle

\addtocounter{footnote}{1}

\section{Introduction}

The integer factorization is an important concept in Mathematics,
 which is related to a lot of  number theory problems, for example, the distribution of primes,  the distribution of almost-primes, the Goldbach conjecture, the distribution of twin primes, etc. Many problems involving the integer factorization have attracted the attention of many authors.

{\bf Definition 1.}
Given an integer $k \geq 2$, we say that a positive integer $n$ is an
{\it almost $k$-th power} if it can be factored as
$$n = n_1n_2\cdots n_k, \ \ \ \  n_1, n_2,\cdots , n_k \asymp n^{1/k}.$$

{\bf Definition 2(see \cite{Ch8}).}      An infinite sequence of positive integers $A$ is called
{\it  almost dense} if, for any
$ \varepsilon > 0$, there exists a constant $c_{\varepsilon,A} >0$ such that
$$  \sum_{\stackrel{X<n\leq  2X }{n\in A}}1 \geq
 c_{\varepsilon,A}X^{1-\varepsilon /4}  $$
for all sufficiently large $ X.$

 In the work \cite{Ch1}-\cite{Ch7},
Chan studied the almost square problem in short intervals and a series of interesting results were obtained.
In 2012,   Islam \cite{Is} showed the existence
of a number $n = n_1n_2 p \in [x, x + x^{1/2}] $ where
$n_1, n_2, p \asymp n^{1/3}, (n_1, n_2, p)=1.$

In 2025,  Chan \cite{Ch8} studied the almost cubes and almost fourth powers
in short intervals. He proved the following Theorem A and Theorem B.

{\bf Theorem A}. For any $\varepsilon>0$ and any two {\it almost dense} sequences $A_1$ and $A_2$, the interval
$[x, x + x^{5/9+\varepsilon}]$ contains an integer $n = m a_1  a_2$ for some
$a_1 \in A_1, a_2 \in A_2$ and integer $m$
with $a_1, a_2,m \asymp x^{1/3}$ for all sufficiently large $x$.

Let $  \mathbb{P} $ denote the set of all primes.
When taking $A_1=A_2=  \mathbb{P},$ Theorem A gives a short interval result  of the distribution  of numbers of the form $p_1p_2m$ with $p_1\asymp p_2\asymp m$.

{\bf Theorem B}. For any $\varepsilon>0$ and any three {\it almost dense} sequences $A_1$, $A_2$ and $A_3$, the interval
$[x, x + x^{34/55+\varepsilon}]$ contains an integer $n = m a_1  a_2a_3$ for some
$a_1 \in A_1, a_2 \in A_2, a_3 \in A_3$ and integer $m$
with $a_1, a_2, a_3, m \asymp x^{1/4}$ for all sufficiently large $x$.

When taking $A_1=A_2= A_3= \mathbb{P},$ Theorem B gives a short interval result  of the distribution  of numbers of the form $p_1p_2p_3m$ with
$p_1\asymp p_2\asymp p_3\asymp m$.

\subsection{\bf Main results}

We first introduce the following definitions.

{\bf Definition 3.}   Suppose $0<\theta\leq 1$ is a fixed real number and $A\subseteq \mathbb{N}$ is an infinite sequence of positive integers.
 We say that $A$ is {\it  $\theta$-short almost dense} if, for any
$ \varepsilon > 0$, there exists a constant $c_{\varepsilon,A} >0$ such that
$$  \sum_{\stackrel{X<n\leq  X+X^\theta}{n\in A}}1 \geq
 c_{\varepsilon,A}X^{1-\varepsilon /4}  $$
for all sufficiently large $ X.$

{\bf Definition $3^{\prime}$.} We can give a more explicit form of the above definition. Suppose $0<\theta\leq 1$ is a fixed real number and $ X$ is any real number   sufficiently large. Suppose $A$ is an infinite sequence of positive integers, whose counting function is denoted by $A(X).$ If there is a function
$F_A: [1, \infty)\rightarrow (0, \infty)$ with $x^{-\varepsilon} \ll F_A(x)\ll x^\varepsilon$ when $x\rightarrow \infty $ such that
the estimate
\begin{equation}
 A(X+X^{\prime})-A(X)= \sum_{\stackrel{X<n\leq  X+X^{\prime}}{n\in A}}1 \geq
   F_A(X)(X^{\prime}-X) \end{equation}
holds for any $X^\theta\ll X^{\prime}\ll X, $ then we say that  $A$ is {\it  $\theta$-short almost dense}.

{\bf Definition 4.} Suppose $0<\theta\leq 1$ is a fixed real number. An integer $n$
is called a $   \theta $-short almost  $k$-th power
 if
 $n = n_1n_2\cdots n_k $ such that $n_j\in [n^{1/k}-C_j n^{\theta/k}, n^{1/k}+C_j n^{\theta/k}] \ (j=1,2, \cdots, k)$ with positive constants
 $C_j>0\ (j=1,2, \cdots, k).$

We are interested in the short interval distribution  of   almost
$k$-th powers. It is interesting to know for how short an interval  must contain an almost $k$-th power. We propose the following two conjectures.

{\bf Conjecture(Weak form).} {\it  Let $k\geq 3$ be a fixed integer and     $A_j\ (1\leq j\leq k-1)$  be    {\it $\theta$-short almost dense} sets satisfying (1.1).
Suppose $0<c_0<1$ is some constant and
$c_0<\theta\leq 1$ is a  fixed real number. Then there exists
a constant $c>0$ such that  in the interval $(x,x+x^{c+ 1/2-c\theta+\varepsilon}] $ there is at least one
  $\theta$-short almost  $k$-th power $n=a_1a_2\cdots a_{k-1}m$   which satisfies
\begin{eqnarray}
 \ \ \ \ a_j\in   A_j, \
   |a_j-n^{1/k}|\ll n^{\theta/k}\ (1\leq j\leq k-1), \ \
    |m-n^{1/k}|\ll n^{\theta/k}.
\end{eqnarray}
    }

{\bf Conjecture(Strong form).} {\it  Let $k\geq 3$ be a fixed integer and     $A_j\ (1\leq j\leq k-1)$  be    {\it $(1-\varepsilon)$-short almost dense} sets satisfying (1.1).
Then there is at least one
  $ (1-\varepsilon) $-short almost  $k$-th power $n=a_1a_2\cdots a_{k-1}m$ in the interval $(x,x+x^{ \varepsilon}] $ such that    (1.2) holds for $\theta=1-\varepsilon$.
   }

Now we   prove the following results,
which confirm that the above weak form conjecture for $k=3$ with $c=1/3$ and $k=4 $ with $c=1/2.$
Note that when $\theta=1,$  Theorem 1 improves
Theorem A and   Theorem 2 improves  Theorem B, respectively.
For any $k\geq 3,$  it seems that the   strong form conjecture is out of reach
of the present methods in Analytic number theory.

{\bf Theorem 1.} {\it  Let   $1/4<\theta\leq 1$ be a fixed real number. Suppose  $A_1$ and $A_2$  are two {\it  $\theta$-short almost dense}   sets satisfying (1.1).
Suppose $x\geq 10$ is a large real number  and
$$  y(\theta)=C\frac{x^{5/6-\theta/3}(\log x)^{3/2}}{F_{A_1}(x^{1/3})F_{A_2}(x^{1/3})},$$
 where $C>0$ is a large positive constant. Then there is at least one
  $ \theta $-short almost  cube $n=a_1a_2m$ in the interval $(x,x+y(\theta)] $ such that
   $$a_j\in A_j,
  |a_j-n^{1/3}|\ll n^{\theta/3}\ (j=1,2),  \  \ |m-n^{1/3}|\ll n^{\theta/3}.$$  }

{\bf Corollary 1.}   {\it Let   $21/40\leq \theta\leq 1$ be a fixed real number.
Suppose $x\geq 10$ is a large real number  and
 $  y(\theta)=C x^{5/6-\theta/3}(\log x)^{7/2},$ where $C>0$ is a large positive constant. Then there is at least one
  $ \theta $-short almost  cube of the form $n=p_1p_2m$ in the interval $(x,x+y(\theta)] $ such that
    $$p_j\in   \mathbb{P},
   |p_j-n^{1/3}|\ll n^{\theta/3}\ (j=1,2), \ \  |m-n^{1/3}|\ll n^{\theta/3}.$$
Especially when $\theta=1,$  there is at least one
 almost  cube of the form $n=p_1p_2m$ in the interval $(x,x+Cx^{1/2}\log^{7/2} x] $ such that $p_1 \asymp p_2\asymp m.$   }

 {\bf Theorem 2.} {\it Suppose $(k-2)/(2k-2)\leq \theta\leq 1$ is a fixed real number, $k\geq 4$ is a fixed integer and     $A_j\ (j=1, 2, \cdots, k-1)$  are {\it $\theta$-short almost dense} sets satisfying (1.1).   Let
\begin{eqnarray}
y(\theta)=x^{\max\left(\frac{3k-4}{2k} -\frac{(2k-4)\theta}{2k},\frac{4k-6}{3k} -\frac{(2k-4)\theta}{3k} \right) }
e^{\kappa c(k-2)\frac{\log x}{\log\log x}}\prod_{j=1}^{k-1}(F_{A_j}(x^{1/k}))^{-1},
\end{eqnarray}
 where $c(\cdot)$ is defined in Lemma 2.5 and $\kappa$ is a fixed constant with
 $(k-1)/k<\kappa<1$.
 Then there is at least one
  $ \theta $-short  $k$-th power of the form $n=a_1a_2\cdots a_{k-1}m$ in the interval $(x,x+y(\theta)] $ such that
    $$a_j\in A_j, \ |a_j-n^{1/k}|\ll n^{\theta/k}\ (j=1,2,\cdots, k-1),  \ \ |m-n^{1/k}|\ll n^{\theta/k}.$$}

{\bf Corollary 2.}   {\it  Let   $21/40\leq \theta\leq 1$ be a fixed real number
and $k\geq 4$ be a fixed integer.
Suppose $x\geq 10$ is a large real number  and let
\begin{eqnarray}
y_k(\theta)=\left\{\begin{array}{ll}
x^{1-\frac{\theta}{2}+\frac{1}{\log\log x}},&\mbox{if $k=4,$}\\
x^{\max\left(\frac{3k-4}{2k} -\frac{(2k-4)\theta}{2k},\frac{4k-6}{3k} -\frac{(2k-4)\theta}{3k} \right)+ \frac{[k/2]  }{\log\log x}}
,& \mbox{if $k\geq 5.$}
\end{array}\right.
\end{eqnarray}
 Then there is at least one
  $ \theta $-short  $k$-th power in the interval $(x,x+y_k(\theta)] $ such that
  $n=p_1p_2\cdots p_{k-1}m$ with
  $$p_j\in   \mathbb{P}, \  \
   |p_j-n^{1/k}|\ll n^{\theta/k}\ (j=1,2,\cdots, k-1), \ \
    |m-n^{1/k}|\ll n^{\theta/k}.$$  }

{\bf Remark 1.} For any $k\geq 3,$ we can compare the short interval distribution
of $k$-th powers and almost $k$-th powers. It is easy  to see that for any $x>1,$ the interval
$(x,x+y]$ contains a $k$-th power if $y>kx^{1-1/k}.$ Note that here the exponent $1-1/k$ is best possible by observing that there are no $k$-th powers in the interval  $(x_0, x_0+y_0]$ with $x_0=[x^{1/k}]^k,\ y_0= ([x^{1/k}]+1/2)^k-[x^{1/k}]^k\asymp x^{1-1/k}$. However, for almost $k$-th powers, we can get much better results. By taking $\theta=1$ in the  above theorems,  we see that for   sufficiently large $x$ and any fixed integer $k\geq 3,$ there exists a constant $\rho_k>0$ such that the interval $(x, x+x^{2/3-\rho_k}]$ contains an almost $k$-th power. Note that
$2/3-\rho_k<1-1/k\ (k\geq 3).$

{\bf Notation.}
Throughout this paper, ${\mathbb  N}$ denotes the set of all positive integers, ${\mathbb  P}$ denotes the set of all primes, ${\mathbb  C}$ denotes the set of all complex numbers. For any $\ell\geq 2,$ $d_{\ell}(n)$ denotes the number of ways $n$ can be written as a product of $\ell$ positive integers. The expression $e(t)$ means  $e(t)=e^{2\pi i t}.$ For any real number
 $t, $ $[t]$ denotes the integer part of  $t,$  $\{t\}$ denotes the fractional part of  $t,$ $\psi(t)=\{t\}-1/2,$ and $\Vert t\Vert=\min(\{t\}, 1-\{t\}).$
 The symbol $n\sim N$ means $N<n\leq 2N$  and $n\asymp N$ means that there exist two absolute positive constants $0<c_1<c_2$ such that $c_1N\leq n\leq c_2N.$ As usual, the symbols $f=O(g)$ and $f\ll g$ mean that $|f|\leq Cg$ for some positive constant $C.$   We always use $\varepsilon$ denote a small positive constant, which maybe different at different places.

\section{Some preliminary lemmas}

In order to prove our results, we need the following lemmas.

{\bf Lemma 2.1.} Let $\mathcal{H}\geq 2$ be any real number. Then
$$\psi(u)= \sum_{1\leq |h|\leq \mathcal{H}}a(h)e(hu) +O\left(\sum_{0\leq |h|\leq \mathcal{H}}b(h)e(hu) \right),$$
where
$$|a(h)|\leq 1/|h| \ (1\leq |h|\leq \mathcal{H}), \ \ \ |b(h)|\leq 1/ \mathcal{H} \ (0\leq |h|\leq \mathcal{H}).$$

\begin{proof}
See Vaaler \cite{V}.
\end{proof}

{\bf Lemma 2.2.} {\it Suppose  $\alpha $ is a non-integer real number. Then
  for any $x\geq 2$ we have
\begin{eqnarray*}
 \sum_{n\leq x   } e(n\alpha)\ll \min\left(x, \frac{1}{\Vert \alpha \Vert}\right).
\end{eqnarray*}}

\begin{proof}
 It is (8.6) of Iwaniec and Kowalski \cite{IK}.
\end{proof}

{\bf Lemma 2.3.} {\it Let $\mathfrak{U}$ and $ \mathfrak{V}$ be two finite sets of real numbers, $\mathfrak{U} \subset[-U,U],
 \mathfrak{V}\subset[-V,V].$ Then for any complex function $\alpha(u)$ and $\beta(v)$ we have
\begin{eqnarray*}
&&|\sum_{u\in  \mathfrak{U}}\sum_{v\in  \mathfrak{V}}\alpha(u)\beta(v)e(uv)|^2 \\
&&\leq 20(1+UV)\sum_{\stackrel {u,u^{'}\in  \mathfrak{U}}{|u-u^{\prime}|\leq V^{-1}}}|\alpha(u)\alpha(u^{'})|
              \sum_{\stackrel {v,v^{'}\in  \mathfrak{V}}{|v-v^{'}|\leq U^{-1}}}|\beta(v)\beta(v^{'})|.
\end{eqnarray*}}

\begin{proof}
  This is Proposition 1 of Fouvry and Iwaniec \cite{FI}.
\end{proof}

{\bf Lemma 2.4.} {\it  Let $M>0,N>0,u_m>0,v_n>0,A_m>0,B_n>0(1\leq m\leq M,1\leq n\leq N),$ and let $Q_1$
and $Q_2$ be given non-negative numbers, $Q_1\leq Q_2$. Then there is a $Q$
such that $Q_1\leq Q\leq Q_2$ and
\begin{eqnarray*}
  \sum_{m=1}^{M}A_{m}Q^{u_m}+\sum_{n=1}^{N}B_{n}Q^{-v_n}\ll \sum_{m=1}^{M}\sum_{n=1}^{N}
(A_{m}^{v_n}B_{n}^{u_m})^{\frac{1}{u_m+v_n}}\\
+\sum_{m=1}^{M}A_{m}Q_{1}^{u_m}+\sum_{n=1}^{N}B_{n}Q_{2}^{-v_n}.
\end{eqnarray*}}

\begin{proof}
  This is Lemma 2.4 of   \cite{GK}.
\end{proof}

{\bf Lemma 2.5.} {\it  Let $k\geq 2$ be a positive integer. Then we have
the estimate
$$d_k(n)\leq e^{c(k)\frac{\log n}{\log\log n}}=
n^{\frac{c(k)}{\log\log n}}\ \ \ (n\geq 10),$$
where
\begin{eqnarray*}
&&c(2)=1, \ c(3)=2 , \ c(k)= [k/2]+1   \ (k\geq 4).
\end{eqnarray*}
}

\begin{proof}
 Lemma 2.5 is essentially formula (1.70) of \cite{I1}. Here we give explicit values of $c(k) $ for $k\geq 2.$

The formula (1.71) of \cite{I1} reads
$$d_2(n)=d(n)\leq e^{C\frac{\log n}{\log\log n}},\ \ \ (n\geq 10),$$
where we $C=(1+\varepsilon/2)\log2.$  For simplicity,
 we take $c(2)=1.$

 When $k=3$ we have
$$d_3(n)=\sum_{\delta|n}d(\delta),$$
which implies that we can take $c(3)=2 .$

When  $k=2\ell $ or  $k=2\ell+1\ (\ell\geq 2)$ we have
$$d_{2\ell}(n)=\sum_{n=n_1n_2\cdots n_{\ell}}d(n_1)d(n_2)\cdots d(n_\ell),$$
$$d_{2\ell+1}(n)=\sum_{n=n_1n_2\cdots n_{\ell}n_{\ell+1}}d(n_1)d(n_2)\cdots d(n_\ell).$$
Since
$$d(n_j)\leq n^{\frac{c(2)}{\log\log n}}\ \ (j=1,2, \cdots, \ell-1), $$
it follows that
 $$d_{2\ell}(n)\leq n^{\frac{(\ell-1) }{\log\log n}} d_3(n)\leq
 n^{\frac{(\ell+1) }{\log\log n}}$$
 and
 $$d_{2\ell+1}(n)\leq n^{\frac{ \ell   }{\log\log n}} d(n)\leq
 n^{\frac{(\ell+1) }{\log\log n}}.$$
 The above two inequalities show that we can take $c(k)=([k/2]+1).$
\end{proof}

{\bf Lemma 2.6.} {\it  Suppose  $X\geq 10$ is a large real number,
      then we have the estimate
$$\sum_{X<p\leq X+X^{21/40}}1\geq \frac{9}{100}\frac{X^{21/40}}{\log X}. $$
}

\begin{proof}
See the proof of Theorem 7.2 in Harman \cite{H}.
\end{proof}

\section{A spacing problem}

In \cite{FI}, Fourvy and Iwaniec studied a spacing problem, which plays an important role in many problems of number theory.
Let $\alpha\beta\not= 0,\  \Delta>0, $ $H, M\geq 1.$ Let
$\mathcal{A}(H,M;\Delta)$ denote the number of the quadruples
$(h_1, h_2, m_1, m_2)$ such that
\begin{equation}
\left|\left(\frac{h_1}{h_2}\right)^\alpha-
\left(\frac{m_1}{m_2}\right)^\beta\right|<\Delta,
\end{equation}
with $H< h_1, h_2\leq 2H, M< m_1, m_2\leq 2M.$ Then
\begin{equation}
\mathcal{A}(H,M;\Delta)\ll HM\log 2HM+\Delta H^2M^2.
\end{equation}
In the proof of (3.2), Fourvy and Iwaniec used the Dirichlet's principle.

In this section, we   generalize  the above spacing problem to the short interval case under the restriction $\alpha=\beta$.

{\bf Lemma 3.1.} {\it Suppose  $ \Theta>0 $  and  $M,N, L_1, L_2$ are  real numbers with $   M\geq 20, 10\leq L_1\leq M, N\geq 1, 1\leq L_2\leq N $.
Let $ \mathcal{A}(M, N, L_1, L_2; \Theta)$ denote    the number of the quadruples
$(m_1, m_2, m_3, m_4)$ such that
\begin{equation}
 \left|m_1- \frac{m_3}{m_4}m_2\right|< \Theta
\end{equation}
with $ M< m_1, m_2\leq M+L_1, N<  m_3, m_4\leq N+L_2.$ Then we have the estimate
\begin{equation}
\mathcal{A}(M, N, L_1, L_2; \Theta) \ll L_1 L_2\log N
 +  N L_2 \log^2 N  +\Theta L_1L_2^2.
 \end{equation}
}

\begin{proof} From (3.3) we have
\begin{equation*}
 \frac{m_3}{m_4}m_2 - \Theta<  m_1\leq  \frac{m_3}{m_4}m_2 + \Theta.
\end{equation*}
Thus we get
\begin{eqnarray}
 \ \ \ \ \ \ \ \ \mathcal{A}(M, N, L_1, L_2; \Theta)
 &&\leq  \sum_{\stackrel{N<  m_3, m_4\leq N+L_2}{M<  m_2\leq M+L_1}}
 \left(\left[\frac{m_3}{m_4}m_2 + \Theta\right]-\left[\frac{m_3}{m_4}m_2- \Theta\right]\right)\\
 &&= \ 2\Theta L_1L_2^2\nonumber\\
 &&\  -\sum_{\stackrel{N < m_3, m_4\leq N+L_2}{M<   m_2\leq M+L_1}}
 \left(\psi\left(\frac{m_3}{m_4}m_2 + \Theta\right)-\psi\left(\frac{m_3}{m_4}m_2- \Theta\right)\right).\nonumber
\end{eqnarray}

So we only need to bound the sum
$$S_{\pm}=\sum_{\stackrel{N<  m_3, m_4\leq N+L_2}{M<   m_2\leq M+L_1}}
  \psi\left(\frac{m_3}{m_4}m_2 \pm \Theta\right).$$

If $N\leq 60,$ or $N_6>60, L_2\leq 30$,
 then trivially we have $S_{\pm}\ll L_1.$
 Later we suppose $N>60, L_2>30.$
Write
\begin{eqnarray}
 S_{\pm}&&=\sum_{1\leq u\leq L_2}
\sum_{\stackrel{\frac{N}{u} < n_3, n_4\leq\frac{N+L_2}{u},(n_3,n_4)=1}{M<  m_2\leq M+L_1}}  \psi\left(\frac{n_3}{n_4}m_2 \pm \Theta\right)\\
&&=S_{\pm, 1}+S_{\pm, 2},\nonumber
\end{eqnarray}
where ($U=L_2/4$)
\begin{eqnarray}
 S_{\pm, 1}&&=\sum_{1\leq u\leq U}
\sum_{\stackrel{\frac{N}{u} < n_3, n_4\leq\frac{N+L_2}{u},(n_3,n_4)=1}{M<  m_2\leq M+L_1}}  \psi\left(\frac{n_3}{n_4}m_2 \pm \Theta\right),\\
 S_{\pm, 2}&&=\sum_{U< u\leq L_2}
\sum_{\stackrel{\frac{N}{u}<  n_3, n_4\leq\frac{N+L_2}{u},(n_3,n_4)=1}{M<   m_2\leq M+L_1}}  \psi\left(\frac{n_3}{n_4}m_2 \pm \Theta\right)\nonumber.
\end{eqnarray}

Trivially we have
\begin{eqnarray}
 S_{\pm, 2} \ll \sum_{U< u\leq L_2}L_1(L_2/u)^2\ll L_1 L_2^2U^{-1}\ll L_1 L_2.
\end{eqnarray}

Let
\begin{eqnarray}
 S(\mathcal{N},\mathcal{L}, L_1)&&=
\sum_{\stackrel{\mathcal{N} < n_3, n_4\leq \mathcal{N}+\mathcal{L},(n_3,n_4)=1}{M<   m_2\leq M+L_1}}  \psi\left(\frac{n_3}{n_4}m_2 \pm \Theta\right), \
  10\leq \mathcal{L}\leq \mathcal{N}. \nonumber
\end{eqnarray}

Suppose $ \mathcal{H}=\mathcal{L}/3.$   By Lemma 2.1 and Lemma 2.2 we
have
\begin{eqnarray}
 \ \ \  S(\mathcal{N},\mathcal{L}, L_1)&&\ll \frac{L_1\mathcal{L}^2}{\mathcal{H}}+\sum_{1\leq h\leq  \mathcal{H}}\frac{1}{h}\left|
  \sum_{\stackrel{\mathcal{N} < n_3, n_4\leq\mathcal{N}+\mathcal{L}, (n_3,n_4)=1}{M<   m_2\leq M+L_1}}
  e\left(\frac{h n_3}{n_4}m_2 + h\Theta\right) \right|\\
 && \ll L_1\mathcal{L} +\sum_{1\leq h\leq  \mathcal{H}}\frac{1}{h}\left|
  \sum_{\stackrel{\mathcal{N} < n_3, n_4\leq\mathcal{N}+\mathcal{L}, (n_3,n_4)=1}{M<  m_2\leq M+L_1}}
  e\left(\frac{h n_3}{n_4}m_2  \right) \right|    \nonumber\\
&&\ll   L_1\mathcal{L}+S,  \nonumber
\end{eqnarray}
where
$$S: =\sum_{1\leq h\leq  \mathcal{H}}\frac{1}{h} \sum_{\mathcal{N}<    n_4<\mathcal{N}+\mathcal{L}}
  \sum_{\stackrel{\mathcal{N}<  n_3\leq\mathcal{N}+\mathcal{L}}{(n_3,n_4)=1 }}
  \frac{1}{\Vert \frac{h n_3}{n_4}   \Vert }.  $$

So we only need to estimate $S.$ Let $d=(h,n_4).$ Write $h=dh_1, n_4=dn_4^{*}, $ then
$n_4^{*}\geq 3, \ (h_1, n_4^{*})=1.$ Thus we have
\begin{eqnarray}
 S =\sum_{d\leq \mathcal{H}}\frac{1}{d} \sum_{1\leq h_1\leq  \frac{\mathcal{H}}{d}}\frac{1}{h_1}
\sum_{ \stackrel{ \frac{\mathcal{N}}{d}<   n_4^{*}\leq \frac{\mathcal{N}}{d}+\frac{\mathcal{L}}{d}}{(h_1, n_4^{*})=1}}
  \sum_{\stackrel{\mathcal{N}<  n_3\leq\mathcal{N}+\mathcal{L}}{(n_3,dn_4^{*})=1 }}
  \frac{1}{\Vert \frac{h_1 n_3}{n_4^{*}}   \Vert }.
   \end{eqnarray}

First we estimate the innermost sum
$$T= \sum_{\stackrel{\mathcal{N}<  n_3\leq\mathcal{N}+\mathcal{L}}{(n_3,dn_4^{*})=1 }}
  \frac{1}{\Vert \frac{h_1 n_3}{n_4^{*}}   \Vert }, $$
  which is bounded by $O(1+\mathcal{L}/n_4^{*})$
  sums of the form
$$T^{*}= \sum_{\stackrel{t<n_3\leq t+n_4^{*} }{(n_3, n_4^{*})=1 }}
  \frac{1}{\Vert \frac{h_1 n_3}{n_4^{*}}   \Vert }, \ t\geq 1.$$

Since $(h_1, n_4^{*})=1,$ we see that when $n_3$ passes through
 a reduced residue system $mod\ n_3$, $h_1n_3$ also  passes through a reduced residue system   $mod\ n_3$. Thus we have
\begin{eqnarray*}
 T^{*}&&= \sum_{\stackrel{1\leq n_3\leq  n_4^{*} }{(n_3, n_4^{*})=1 }}
  \frac{1}{\Vert \frac{  n_3}{n_4^{*}}   \Vert }
 \ll \sum_{\stackrel{1\leq n_3\leq  n_4^{*}-1 }{(n_3, n_4^{*})=1 }}
  \frac{  n_4^{*}}{     n_3  } \ll n_4^{*}\log n_4^{*}.
\end{eqnarray*}

So we get
\begin{eqnarray*}
 T\ll  (1+\mathcal{L}/n_4^{*})\times n_4^{*}\log n_4^{*}\ll
 ( \mathcal{L}+n_4^{*}) \log n_4^{*}.
\end{eqnarray*}

Inserting the above bound into (3.10) we get
\begin{eqnarray}
 S&&\ll \sum_{d\leq \mathcal{H}}\frac{1}{d} \sum_{1\leq h_1\leq  \frac{\mathcal{H}}{d}}\frac{1}{h_1}
\sum_{ \stackrel{ \frac{\mathcal{N}}{d}<   n_4^{*}\leq \frac{\mathcal{N}}{d}+\frac{\mathcal{L}}{d}}{(h_1, n_4^{*})=1}}
( \mathcal{L}+n_4^{*}) \log n_4^{*}\\
&&\ll \sum_{d\leq \mathcal{H}}\frac{1}{d} \sum_{1\leq h_1\leq  \frac{\mathcal{H}}{d}}\frac{1}{h_1}
\left( \frac{\mathcal{L}^2}{d}+\frac{\mathcal{N}\mathcal{L} }{d^2}\right) \log \mathcal{N}\nonumber\\
&&\ll  (\mathcal{L}^2 + \mathcal{N}\mathcal{L}  ) \log^2 \mathcal{N}
\ll    \mathcal{N}\mathcal{L}   \log^2 \mathcal{N}.\nonumber
\end{eqnarray}

From (3.9) and (3.11) we get
\begin{eqnarray}
S(\mathcal{N},\mathcal{L}, L_1) \ll     L_1\mathcal{L}+\mathcal{N}\mathcal{L}
  \log^2 \mathcal{N}.
\end{eqnarray}

From (3.7) and (3.12) we get
\begin{eqnarray}
 S_{\pm, 1} \ll \sum_{1\leq u\leq U}\left(L_1\frac{ L_2}{u}+\frac{N L_2}{u^2}\log^2 \mathcal{N} \right)\ll L_1 L_2\log N + N L_2 \log^2 N.
\end{eqnarray}

From (3.5), (3.6), (3.8) and (3.13) we get the estimate (3.4).
This completes the proof of Lemma 3.1.
\end{proof}

{\bf Lemma 3.2.} {\it  Let $ \beta\eta\not= 0$ be   fixed real numbers,
$0< \eta\leq 1, $ $\  \Delta>0,$  and  $ H\geq 1,  M\geq 20, M\geq H$. Let
$\mathcal{A}(H,M;\Delta; \eta)$ denote the number of the quadruples
$(h_1, h_2, m_1, m_2)$ such that
\begin{equation}
\left|\left(\frac{h_1}{h_2}\right)^\beta-
\left(\frac{m_1}{m_2}\right)^\beta\right|<\Delta,
\end{equation}
with $H<h_1, h_2\leq 2H, M <m_1, m_2\leq M+M^\eta.$ Then
\begin{equation}
\mathcal{A}(H,M;\Delta; \eta)
\ll H M^\eta \log 2H  +H^2 \log^2 2H  +\Delta M^{1+\eta} H^2.
\end{equation} }

\begin{proof} Using Lagrange theorem   we have
\begin{eqnarray*}
  \left(\frac{m_1}{m_2}\right)^\beta  - \left(\frac{h_1}{h_2}\right)^\beta
 =\beta u_0^{\beta-1}\left(\frac{m_1}{m_2}- \frac{h_1}{h_2}  \right)=
\frac{\beta u_0^{\beta-1}}{m_2}\left( m_1- \frac{h_1}{h_2}m_2\right),
\end{eqnarray*}
where
$1/2<u_0<2.$ Hence $ \mathcal{A}(H, M; \Delta, \eta)$ doesn't exceed the number of  the quadruples $(m_1, m_2, h_1, h_2)$ such that
$$ \left| m_1- \frac{h_1}{h_2}m_2\right|< C(\beta) \Delta M,\ \ H<h_1, h_2\leq 2H, M <m_1, m_2\leq M+M^\eta, $$
where $C(\beta)$ is a positive constant.
Now Lemma 3.2  follows from
Lemma 3.1 by taking
$(\Theta, M, N, L_1, L_2)=(C(\beta) \Delta M, M,  H,  M^\eta, H)$.
\end{proof}

{\bf Remark 3.1}  If  $\alpha=\beta,$ then (3.1) becomes (3.14).
Take $\eta=1 $ in Lemma 3.2. If $H\ll M\log^{-1} 2M, $
 the estimate (3.15) becomes
\begin{eqnarray*}
\mathcal{A}(H,M;\Delta; 1)&&
\ll H M  \log 2H  +H^2 \log^2 2H  +\Delta M^2 H^2\\
&&\ll H M  \log 2HM     +\Delta M^2 H^2.
\end{eqnarray*}
  So Lemma 3.2 provides a different proof of (3.2) under the condition
  $\alpha=\beta$ and  $H\ll M\log^{-1} 2M  $ (the latter is often not important in applications).

\section{Estimate  of an exponential sum}

The exponential sum  of the form
\begin{eqnarray*}
 \sum_{h\sim H}\sum_{m\sim m }
\sum_{  n\sim N }\rho_h a_{m} b_{n}
e\left(\frac{hX}{mn }\right) \ \ (\rho_h, a_{m}, b_{n}\in {\Bbb C})
\end{eqnarray*}
has many applications in number theory.
In this section, we shall estimate  the exponential sum  of the above form with
some variables being restricted in short intervals.

Suppose $x$ and $M_1, M_2$ are large real numbers with $x^{\varepsilon}\ll
 M_1,M_2, M_1M_2\ll x^{1-\varepsilon}, $  $0<\eta_1\leq 1, 0<\eta_2\leq 1$ are fixed real numbers and suppose $1\leq H\leq  M_1$.
Define
\begin{eqnarray*}
&&S(x, H, M_1,M_2; \eta_1, \eta_2):=\sum_{h\sim H}\rho_h
\sum_{ M_1 < m_1\leq M_1+M_1^{\eta_1} }a_{m_1}
\sum_{ M_2 < m_2\leq M_2+M_2^{\eta_2} }b_{m_2}
e\left(\frac{hx}{m_1m_2 }\right),
\end{eqnarray*}
where $\rho_h, a_{m_1}, b_{m_2}\in {\Bbb C}$ with $\rho_h\ll 1, a_{m_1}\ll 1, b_{m_2}\ll 1.$ Let $F : =Hx/M_1M_2$ and $  \mathfrak{L}=\log x.$  Then we have the following proposition.

{\bf Proposition 4.1.} {\it Under the above notation and suppose $F \gg M_2$, we have the estimate
\begin{eqnarray*}
 S(x, H, M_1,M_2; \eta_1, \eta_2)&&
\ll   H^{1/2}F^{1/2} M_1^{\eta_1/2}M_2^{\eta_2/2} \mathfrak{L}^{1/2} \\
  &&\ \ \  +F^{1/2}H M_2^{\eta_2/2}\mathfrak{L} +  H^{3/2} M_1^{\eta_1/2}M_2^{\eta_2/2}.\nonumber
\end{eqnarray*}}

\begin{proof}
Let
\begin{eqnarray*}
&&U=4H/M_1, \ \ V=x/M_2,\\
&&\mathfrak{U}=\{h/m_1: H< h\leq 2H, M_1< m_1\leq M_1+M_1^{\eta_1}\},\\
&&\mathfrak{V}=\{x/m_2:   M_2< m_2\leq M_2+M_2^{\eta_2}\}.
\end{eqnarray*}

By Lemma 2.3 we have
\begin{eqnarray}
&&|S(x, H, M_1,M_2; \eta_1, \eta_2)|^2 \ll F\mathcal{B}_1\mathcal{B}_2,
\end{eqnarray}
where $\mathcal{B}_1$ stands for the number of the quadruples
$(h_1, h_2, m_1, \widetilde{m}_1)$ such that
\begin{equation*}
\left| \frac{h_1}{h_2} - \frac{m_1}{\widetilde{m}_1} \right|<\frac{M_2}{x},\
H< h_1, h_2\leq 2H,\ M_1< m_1, \widetilde{m}_1\leq M_1+M_1^{\eta_1}
\end{equation*}
and
 $\mathcal{B}_2$ stands for the number of pairs
$(m_2, \widetilde{m}_2)$ such that
\begin{equation*}
\left| \frac{1}{m_2} - \frac{1}{\widetilde{m}_2} \right|<\frac{M_1}{Hx},\ \
M_2< m_2, \widetilde{m}_2\leq M_2+M_2^{\eta_2}.
\end{equation*}

Taking $(H, M, \Delta, \eta)=(H, M_1, M_2/x, \eta_1)$ in Lemma 3.2 we get
\begin{eqnarray}
 \mathcal{B}_1&&
\ll H M_1^{\eta_1} \log 2H  +H^2 \log^2 2H  +M_2x^{-1} M_1^{1+\eta_1} H^2\\
&&\ll H M_1^{\eta_1}  \mathfrak{L}   +H^2\mathfrak{L}^2 + F^{-1}H^3 M_1^{\eta_1}.\nonumber
\end{eqnarray}

It is easy to see that $\mathcal{B}_2$ does not exceed the number of solutions of the
inequality
\begin{equation*}
 |  m_2 -  \widetilde{m}_2   |<\frac{4M_1M_2^2}{Hx}=\frac{4 M_2 }{F}.
\end{equation*}
So we have
\begin{eqnarray}
\mathcal{B}_2 \ll M_2^{\eta_2}+ \frac{M_1M_2^{2+\eta_2}}{Hx}
= M_2^{\eta_2}+ \frac{ M_2^{1+\eta_2}}{F}
\ll M_2^{\eta_2}
\end{eqnarray}
by noting $F\gg M_2.$

Now Proposition 4.1 follows from (4.1)-(4.3).
\end{proof}


\section{Proofs of results}

\subsection{Preparation}

Let $k\geq 3$ be a fixed positive integer.
Suppose $0<\theta\leq 1$ is a fixed real number and
suppose $A_j\subseteq {\Bbb N}$ ($j=1, 2, \cdots, k-1$)
are {\it  $\theta$-short almost dense} sets.
Let $0<c\leq 1/2$ be a fixed small constant. Define
\begin{equation*}
 f_k(n):=f(n;A_1,\cdots, A_{k-1}; \theta, c)= \sum_{\substack{n=n_1n_2\cdots n_{k-1}
  n_k\\n_j\in I_k(n;\theta,c)\ (1\leq j\leq k-1)\\n_j\in A_j\ (1\leq j\leq k-1)}}1,
\end{equation*}
 where $
I_k(n;\theta,c):=[n^{1/k}-cn^{\theta/k}, n^{1/k}+cn^{\theta/k}].$

 We first show that if $f_k(n)\geq 1,$ then $n$ must be a {\it $\theta$-short almost} $k$-th power with $n_j\in A_j\ (1\leq j\leq k-1).$  If $\theta=1,$ then $I_k(n;1,c)=[(1-c)n^{1/k}, (1+c)n^{1/k} ]$ and
 we have
 $$n_k=\frac{n}{n_1n_2\cdots n_{k-1}}\in \left[(1+c)^{-(k-1)}n^{1/k}, (1-c)^{-(k-1)}n^{1/k}\right].$$
So $n$ is a $1$-short {\it almost} $k$-th power. Now suppose $0<\theta< 1 $ and in this case we have
 $$ \frac{n}{(n^{1/k}+cn^{\theta/k})^{k-1}}\leq n_k\leq \frac{n}{(n^{1/k}-cn^{\theta/k})^{k-1}}.
 $$
 It is easy to see that
 \begin{eqnarray*}
 && \frac{n}{(n^{1/k}+cn^{\theta/k})^{k-1}}
  =n^{1/k}- c(k-1)n^{\theta/k}+ O(n^{\theta/k-(1-\theta)/k}),\\
  && \frac{n}{(n^{1/k}-cn^{\theta/k})^{k-1}}
  =n^{1/k}+c(k-1)n^{\theta/k}+ O(n^{\theta/k-(1-\theta)/k}).
 \end{eqnarray*}
So it follows that $n$ is a $\theta$-short {\it almost} $k$-th power.

Suppose $x>y>1$ are large real numbers such that
\begin{eqnarray}
  yx^{ 1/k-1}\ll x^{\theta/k}
\end{eqnarray}
and define
$$S_k(x):= \sum_{x<n\leq x+y}f_k(n). $$

Let
$$X_1^{(k)}=x^{1/k}-c x^{\theta/k}, \ X_2^{(k)}= (x+y)^{1/k}+c (x+y)^{\theta/k}. $$
By (5.1) and Lagrange's theorem we have
\begin{eqnarray}
x^{\theta/k}\ll X_2^{(k)}-X_1^{(k)} =  (x+y)^{1/k}-x^{1/k}+c (x+y)^{\theta/k}+c x^{\theta/k}
\ll x^{\theta/k}.
\end{eqnarray}

We can write
\begin{eqnarray}
S_k(x)&& = \sum_{\stackrel{X_1^{(k)}<n_1\leq X_2^{(k)}}{n_1\in A_1 }}\cdots
\sum_{\stackrel{X_1^{(k)}<n_{k-1}\leq X_2^{(k)}}{n_{k-1}\in A_{k-1} }}
\sum_{\frac{x}{n_1\cdots n_{k-1}}<n_k\leq \frac{x+y}{n_1\cdots n_{k-1}}}1\\
&&=\sum_{\stackrel{X_1^{(k)}<n_1\leq X_2^{(k)}}{n_1\in A_1 }}\cdots
\sum_{\stackrel{X_1^{(k)}<n_{k-1}\leq X_2^{(k)}}{n_{k-1}\in A_{k-1} }}
\left(\left[\frac{x+y}{n_1\cdots n_{k-1}}\right]-\left[\frac{x}{n_1\cdots n_{k-1}}\right]\right)\nonumber\\
&&=M_k(x;y)+E_k(x;y),\nonumber
\end{eqnarray}
where
\begin{eqnarray*}
M_k(x;y) &&=y\prod_{j=1}^{k-1} \sum_{\stackrel{X_1^{(k)}<n_j\leq X_2^{(k)}}{n_1\in A_j }}\frac{1}{n_j},\\
E_k(x;y)&&=\sum_{\stackrel{X_1^{(k)}<n_1\leq X_2^{(k)}}{n_1\in A_1 }}\cdots
\sum_{\stackrel{X_1^{(k)}<n_{k-1}\leq X_2^{(k)}}{n_{k-1}\in A_{k-1} }}
\left(\psi\left(\frac{x+y}{n_1\cdots n_{k-1}}\right)-\psi\left(\frac{x}{n_1\cdots n_{k-1}}\right)\right)\\
&&=\sum_{j=0}^1\sum_{\stackrel{X_1^{(k)}<n_1\leq X_2^{(k)}}{n_1\in A_1 }}\cdots
\sum_{\stackrel{X_1^{(k)}<n_{k-1}\leq X_2^{(k)}}{n_{k-1}\in A_{k-1} }}(-1)^{j }
 \psi\left(\frac{x+\chi(j)y}{n_1\cdots n_{k-1}}\right),
\end{eqnarray*}
where $\chi(\cdot)$ is the characteristic function mod $2.$

Recall that $A_j\subseteq {\Bbb N}$ ($j=1, 2, \cdots, k-1$)
are {\it  $\theta$-short almost dense} sets. According to (1.1), we have
\begin{eqnarray*}
\sum_{\stackrel{X_1^{(k)}<n_j\leq X_2^{(k)}}{n_j\in A_j }} 1\gg F_j(X_1^{(k)})x^{\theta/k}, \ (j=1,2, \cdots, k-1)
\end{eqnarray*}
with
$\ \
X^{-\varepsilon} \ll F_j(X)\ll X^\varepsilon, \ \ X\rightarrow \infty.
$

Thus we get
\begin{eqnarray}
&&M_k(x;y)\gg y x^{-\frac{k-1}{k}(1-\theta)} \prod_{j=1}^{k-1}F_j(X_1^{(k)}) .
\end{eqnarray}

So now we only need to bound the sum $E_k(x;y).$ Let $$M=X_1^{(k)}\asymp x^{1/k},  L=X_2^{(k)}-X_1^{(k)}. $$  From  (5.2) we have
\begin{eqnarray}
 x^{\frac{\theta}{k}}\ll L\ll  x^{\frac{\theta}{k}}.
\end{eqnarray}
Suppose $\mathcal{H}_k\geq 10$ is a real parameter to be determined later. By Lemma 2.1 and (5.5) we have
\begin{eqnarray}
&&\ \ \ \ E_k(x;y)\\&&\ll \frac{L^{k-1}}{\mathcal{H}_k}
 +\sum_{j=0}^1\sum_{1\leq |h|\leq \mathcal{H}_k }\frac{1}{|h|}\left|\sum_{\stackrel{X_1^{(k)}<n_1\leq X_2^{(k)}}{n_1\in A_1 }}\cdots
\sum_{\stackrel{X_1^{(k)}<n_{k-1}\leq X_2^{(k)}}{n_{k-1}\in A_{k-1} }}
 e\left(\frac{hx+h\chi(j)y}{n_1\cdots n_{k-1}}\right)\right|\nonumber\\
 &&\ll \frac{L^{k-1}}{\mathcal{H}_k}
 + \sum_{1\leq  h \leq \mathcal{H}_k }\frac{1}{ h }\left|\sum_{\stackrel{X_1^{(k)}<n_1\leq X_2^{(k)}}{n_1\in A_1 }}\cdots
\sum_{\stackrel{X_1^{(k)}<n_{k-1}\leq X_2^{(k)}}{n_{k-1}\in A_{k-1} }}
 e\left(\frac{hx }{n_1\cdots n_{k-1}}\right)\right|\nonumber\\
 &&\ll \frac{L^{k-1}}{\mathcal{H}_k}
 + \frac{\log\mathcal{H}_k}{H}\times S_{k}(H,x,\theta,\mathfrak{A}_k)   \nonumber
\end{eqnarray}
for some $1\ll H\ll \mathcal{H}_k,$ where
$$S_{k}(H,x,\theta,\mathfrak{A}_k):=\sum_{   h \sim H }\left|\sum_{\stackrel{X_1^{(k)}<n_1\leq X_2^{(k)}}{n_1\in A_1 }}\cdots
\sum_{\stackrel{X_1^{(k)}<n_{k-1}\leq X_2^{(k)}}{n_{k-1}\in A_{k-1} }}
 e\left(\frac{hx }{n_1\cdots n_{k-1}}\right)\right|,$$
with $\mathfrak{A}_k=(A_1,\cdots, A_{k-1}).$

\subsection{Proof of Theorem 1 }\

When $k=3$ we have
 \begin{eqnarray*}
S_{3}(H,x,\theta,\mathfrak{A}_3)= \sum_{   h \sim H } \rho(h)\sum_{\stackrel{X_1^{(3)}<n_1\leq X_2^{(3)}}{n_1\in A_1 }}
\sum_{\stackrel{X_1^{(3)}<n_{2}\leq X_2^{(3)}}{n_{2}\in A_{2} }}
 e\left(\frac{hx }{n_1  n_{2}}\right),
\end{eqnarray*}
where $\rho(h)\ll 1.$

Take $M_1=M_2=M=X_1^{(3)}\asymp x^{1/3}, \eta_1=\eta_2=\theta, F=Hx/M_1M_2\asymp Hx^{1/3}\gg M_2 $ in Proposition 4.1. We get
\begin{eqnarray}\ \ \ \ \
S_{3}(H,x,\theta,\mathfrak{A}_3)\ll H x^{1/6+\theta/3}  (\log x)^{1/2}
     +x^{1/6+\theta/6}H^{3/2} \log x +      H^{3/2} x^{\theta/3}.
\end{eqnarray}

Inserting the bound (5.7) into (5.6) with $k=3,$ we get
\begin{eqnarray}
\ \ \ \ \ \ \ \ \  E_3(x;y)  \ll   \frac{L^2}{\mathcal{H}_3}
 +    x^{1/6+\theta/3}  (\log x)^{3/2}
     +\mathcal{H}_3^{1/2}x^{1/6+\theta/6} \log^2 x +      \mathcal{H}_3^{1/2} x^{\theta/3}\log x.
\end{eqnarray}
The   formula (5.8) obviously holds for $0<\mathcal{H}_3\leq 10.$ So taking a best
$\mathcal{H}_3\in (0, x^{1/3})$ by Lemma 2.4 we get
\begin{eqnarray}
   E_3(x;y) && \ll       x^{1/6+\theta/3}  (\log x)^{3/2}
     + x^{1/9+\theta/3} \log^{4/3}x +        x^{4\theta/9}(\log x)^{2/3}\\
     &&\ll         x^{1/6+\theta/3}  (\log x)^{3/2}.\nonumber
\end{eqnarray}

When $k=3$, the formula (5.4) reads
\begin{eqnarray}
&&M_3(x;y)\gg y x^{-\frac{2}{3}(1-\theta)}  F_1(x^{1/3}) F_2(x^{1/3}).
\end{eqnarray}

From (5.1) with $k=3$, (5.9) and (5.10) we get the proof of Theorem 1 for the case $k=3.$

\subsection{Proof of Theorem 2 }\

When $k\geq 4$ we have
$$S_{k}(H,x,\theta,\mathfrak{A}_k)=\sum_{h \sim H } \rho(h)\sum_{\stackrel{X_1^{(k)}<n_1\leq X_2^{(k)}}{n_1\in A_1 }}\cdots
\sum_{\stackrel{X_1^{(k)}<n_{k-1}\leq X_2^{(k)}}{n_{k-1}\in A_{k-1} }}
 e\left(\frac{hx }{n_1\cdots n_{k-1}}\right),$$
where $ \rho(h)\ll 1.$

For  any set $B\subseteq \mathbb{N},$ let
\begin{eqnarray*}
 \xi_{B}(n)=\left\{\begin{array}{ll}
1,&\mbox{if $X_1^{(k)}<n\leq X_2^{(k)}$ and $n\in B$},\\
0 ,& \mbox{if otherwise.}
\end{array}\right.
\end{eqnarray*}

 Write
\begin{eqnarray*}
&& m_1=n_1\cdots n_{k-2}, \ X_1^{(k)}<n_j\leq X_2^{(k)},\ n_j\in A_j\ (1\leq j\leq k-2),\\
&&a_{m_1}=\sum_{  m_1=n_1\cdots n_{k-2}}\xi_{A_1}(n_1) \cdots \xi_{A_{k-2}}(n_{k-2}).
\end{eqnarray*}

If  $X_1^{(k)}<n_j\leq X_2^{(k)}\ (1\leq j\leq k-2),$ then
we have
$$x^{\frac{k-2}{k}}\ll (X_1^{(k)})^{k-2}<m_1
\leq (X_2^{(k)})^{k-2}\ll (X_1^{(k)})^{k-2}+x^{1+\frac{\theta}{k}-\frac{3}{k}}$$
and
$$ x^{1+\frac{\theta}{k}-\frac{3}{k}} \asymp \left((X_1^{(k)})^{k-2}\right)^{\frac{k+\theta-3}{k-2}}.$$

By Lemma 2.5 we have
$$a_{m_1}\leq d_{k-2}(m_1)\ll e^{c(k-2)\frac{\log x^{(k-2)/k}}{\log\log x^{(k-2)/k}}}\ll
e^{c(k-2)\frac{\log x }{\log\log x}},$$
where $c(\cdot)$ was defined in Lemma 2.5.

So we get that
\begin{eqnarray*}
  S_{k}(H,x,\theta,\mathfrak{A}_k)\ll
e^{\frac{(k-2)c(k-2)}{k}\frac{\log x}{\log\log x}}|S^{*}(x, H, M_1,M_2; \eta_1, \eta_2)|,
\end{eqnarray*}
where
\begin{eqnarray*}
&&\ \ \ \ S^{*}(x, H, M_1,M_2; \eta_1, \eta_2)\\
&&=
\sum_{h\sim H}\rho_h
\sum_{ M_1<  m_1\leq M_1+M_1^{\eta_1} }a_{m_1}^{*}
\sum_{ M_2<  m_2\leq M_2+M_2^{\eta_2} }\xi_{A_{k-1}}(m_2)
e\left(\frac{hx}{m_1m_2 }\right)
\end{eqnarray*}
with
\begin{eqnarray*}
&&M_1=(X_1^{(k)})^{k-2}, \ M_2=X_1^{(k)},\ \eta_1=\frac{k+\theta-3}{k-2},
 \\&& \eta_2=\theta,\   F_1\gg  M_2, \ a_{m_1}^{*}=a_{m_1}
 e^{-\frac{(k-2)c(k-2)}{k}\times \frac{\log x}{\log\log x}}\ll 1.
\end{eqnarray*}
We apply  Proposition 4.1 to bound the sum
$S^{*}(x, H, M_1,M_2; \eta_1, \eta_2)$ and get
\begin{eqnarray*}
S_{k}(H,x,\theta,\mathfrak{A}_k)e^{-\frac{(k-2)c(k-2)}{k}\frac{\log x}{\log\log x}} &&\ll
 H^{1/2}F^{1/2} M_1^{\eta_1/2}M_2^{\eta_2/2} \log^{1/2} x \\
  &&\ \ \  +F^{1/2}H M_2^{\eta_2/2}\log  x  +  H^{3/2} M_1^{\eta_1/2}M_2^{\eta_2/2}\\
 &&\ll Hx^{\frac{k-2+2\theta}{2k}}\log^{1/2} x+H^{3/2}x^{\frac{1+\theta}{2k}} \log x
 +H^{3/2}x^{\frac{k-3+2\theta}{2k}}\\
 &&\ll Hx^{\frac{k-2+2\theta}{2k}}\log^{1/2} x
 +H^{3/2}x^{\frac{k-3+2\theta}{2k}}.
\end{eqnarray*}

Inserting the above bound   into (5.6) and recalling (5.5)   we get
\begin{eqnarray*}
  E_k(x;y)e^{- \frac{(k-1)c(k-2)}{k} \frac{\log x}{\log\log x}}\ll   \frac{x^{\frac{(k-1)\theta}{k}}}{\mathcal{H}_k}+ x^{\frac{k-2+2\theta}{2k} }+
\mathcal{H}_k^{1/2}x^{\frac{k-3+2\theta}{2k} }.
\end{eqnarray*}

It is obvious that
the above formula holds for $0<\mathcal{H}_k\leq 10.$ So taking a best
$\mathcal{H}_k \in  (0, x^{(k-2)/k})$ by Lemma 2.4 we get
\begin{eqnarray}
   E_k(x;y) \ll       \left( x^{\frac{k-2+2\theta}{2k} }+
 x^{\frac{ k-3+(k+1)\theta}{3k} }\right) e^{\frac{ (k-1)c(k-2)}{k}\frac{\log x}{\log\log x}}.
\end{eqnarray}

Now  Theorem 2 follows from (5.1), (5.4) and (5.11).

\subsection{Proofs of Corollaries}\

From Lemma 2.6 we see that if
  $X^{21/40}\leq X^{\prime}\leq X,$ then
\begin{equation}
 \mathbb{P}(X+X^{\prime})-\mathbb{P}(X)=
 \sum_{\stackrel{X<p\leq  X+X^{\prime}}{p\in \mathbb{P}}}1 \gg
    F_\mathbb{P}(X)(X^{\prime}-X)
\end{equation}
with $F_\mathbb{P}(X)=\log^{-1} X.$
So the set of all primes $\mathbb{P}$ is a $21/40$-short almost dense set.

If $k=3$,  we get Corollary 1  easily  from Theorem 1 and  (5.12) .

If $k\geq 4,$ the value of $y(\theta)$ in Theorem 2 is
\begin{eqnarray*}
y(\theta)&&=x^{\max\left(\frac{3k-4}{2k} -\frac{(2k-4)\theta}{2k},\frac{4k-6}{3k} -\frac{(2k-4)\theta}{3k} \right) }
e^{\kappa c(k-2)\frac{\log x}{\log\log x}}\prod_{j=1}^{k-1}(F_{A_j}(x^{1/k}))^{-1}\\
&&\asymp x^{\max\left(\frac{3k-4}{2k} -\frac{(2k-4)\theta}{2k},\frac{4k-6}{3k} -\frac{(2k-4)\theta}{3k} \right) }
e^{\kappa c(k-2)\frac{\log x}{\log\log x}} (\log x)^{k-1}\\
&&=x^{\max\left(\frac{3k-4}{2k} -\frac{(2k-4)\theta}{2k},\frac{4k-6}{3k} -\frac{(2k-4)\theta}{3k} \right) }
e^{\kappa c(k-2)\frac{\log x}{\log\log x}} (\log x)^{k-1}\\
&&=x^{\max\left(\frac{3k-4}{2k} -\frac{(2k-4)\theta}{2k},\frac{4k-6}{3k} -\frac{(2k-4)\theta}{3k} \right) }
e^{  c(k-2)\frac{\log x}{\log\log x}}\times \frac{(\log x)^{k-1}}{
e^{ (1-\kappa) c(k-2)\frac{\log x}{\log\log x}}}\\
&&\ll x^{\max\left(\frac{3k-4}{2k} -\frac{(2k-4)\theta}{2k},\frac{4k-6}{3k} -\frac{(2k-4)\theta}{3k} \right) }
e^{  c(k-2)\frac{\log x}{\log\log x}}.
\end{eqnarray*}

So we get Corollary 2 with the help of the values $c(\cdot)$ in Lemma 2.5.

\bigskip

{\bf Conflict of Interest.} The author declares that there is no conflict of interests
regarding the publication of this paper.

\end{document}